\documentclass[12pt,leqno]{amsart}
\usepackage{amsmath, amsfonts}
\usepackage{graphicx,color}
\usepackage{hyperref}
\usepackage[latin1]{inputenc}

\newcommand{\R}{\mathbb{R}}
\def\C{\mathbb{C}}
\def\N{\mathbb{N}}

\theoremstyle{plain}
\newtheorem{theo}{Theorem}[section]

\newtheorem{lemma}[theo]{Lemma}
\newtheorem{prop}[theo]{Proposition}

\theoremstyle{definition}
\newtheorem{defi}[theo]{Definition}

\newtheorem{rem}[theo]{Remark}

\numberwithin{equation}{section}

\def\a{\alpha}

\def\L{\mathbb{L}}
\def\Log{\mathop{\rm Log}\nolimits}
\def\bM{\mathbb{M}}
\def\M{\mathbb{M}}
\newcommand{\hM}{\widehat{\M}}

\def\Re{\mathop{\rm Re}\nolimits}      
\def\Im{\mathop{\rm Im}\nolimits}

\newenvironment{proof1}{\medskip\par\noindent{\bf Proof}\ \normalfont}{\hfill $\Box$
\medskip\par}

\begin{document}

\title[Extension operators for some ultraholomorphic classes of rapid growth]{Extension operators for some ultraholomorphic classes defined by sequences of rapid growth}

\author{Javier Jim\'enez-Garrido\and Alberto Lastra\and Javier Sanz}

\date{\today}

\maketitle
\begin{abstract}
While the asymptotic Borel mapping, sending a function into its series of asymptotic expansion in a sector, is known to be surjective for arbitrary openings in the framework of ultraholomorphic classes associated with sequences of rapid growth, there is no general procedure to construct extension operators in this case. We do provide such operators in complex sectors for some particular classes considered by S.~Pilipovi{\'c}, N.~Teofanov and F.~Tomi{\'c} in the ultradifferentiable setting. Although these classes are, in their words, ``beyond Gevrey regularity'', in some cases they keep the property of stability under differentiation, which is crucial for our technique, based on formal Borel- and truncated Laplace-like transforms with suitable kernels.
\end{abstract}

\section{Introduction}

The concept of asymptotic expansion dates back to H. Poincar\'e in 1886, who gave an analytical meaning to the formal (generally divergent) power series solutions to some classical problems in optics or astronomy. Since then, its use has pervaded the theory of formal solutions to general classes of ordinary and partial differential equations, difference and $q$-difference equations, and so on.
For simplicity, we consider asymptotic expansions $\sum_{p=0}^\infty a_pz^p$ of functions $f$ analytic in unbounded sectors $S$ of the Riemann surface of the logarithm with vertex at 0. Although in its original formulation the remainders $f(z)-\sum_{p=0}^{n-1}a_pz^p$ were assumed to be just bounded for all $z\in S$ by a suitable constant (depending on $n$) times $|z|^n$, the asymptotic techniques have shown its real power when a precise rate of growth for the constants, in the form $CA^nM_n$ for $A,C>0$ and a sequence $(M_p)_{p=0}^\infty$ of positive real numbers, is either obtained from the concrete problem under study, or specified for theoretical considerations. In this way, ultraholomorphic classes of functions associated with the weight sequence $\M=(M_p)_{p=0}^\infty$ and the sector $S$ have been introduced, and it is important to decide about the injectivity and surjectivity of the (asymptotic) Borel map,  sending a function in one of such classes into its series of asymptotic expansion. We refer the reader to the works~\cite{Mandelbrojt,Salinas,JimenezSanzSchindlInjectSurject} for
a complete solution of the injectivity problem. In the last of these references, the seminal results of V. Thilliez~\cite{Thilliez03}, proving surjectivity of the Borel map in a constructive way for strongly regular sequences $\M$ in sectors of opening less than $\pi\gamma(\M)$, where $\gamma(\M)$ is a growth index associated with $\M$, were shown to be optimal except possibly for the critical opening: if surjectivity occurs, the opening will be at most $\pi\gamma(\M)$ (see Subsection~\ref{subsec.WeightSeq} for the notions related to weight sequences). Very recently, G. Schindl and the first and third authors~\cite{JimenezSanzSchindlSurjectDc} have obtained a similar result
for the more general case that $\hM=(p!M_p)_{p=0}^\infty$ is a regular weight sequence in the sense of E.~M.~Dyn'kin~\cite{Dynkin80}, but now the proof is not constructive. Indeed, the key fact is due to A. Debrouwere~\cite{momentsdebrouwere}, who characterized the surjectivity of the Stieltjes moment problem in Gelfand-Shilov spaces defined by regular sequences by abstract functional-analytic techniques. The Fourier transform allows to transfer this information into the asymptotic framework in a halfplane, and in~\cite{JimenezSanzSchindlSurjectDc} Laplace and Borel transforms of arbitrary order extend the procedure to general sectors. In particular, the condition $\gamma(\M)=\infty$ (we say the sequences satisfying this condition are of rapid growth, see Remark~\ref{rem.RapidGrowth}) is shown to amount to the surjectivity of the Borel map, and the existence of both local and global linear continuous right inverses for it, for unbounded sectors of arbitrary (finite) opening in the Riemann surface of the logarithm.

The main aim of this paper is to give a constructive proof of the surjectivity of the Borel map in sectors of the complex plane for the ultraholomorphic class associated with some specific sequences, $\M^{\tau,\sigma}=(p^{\tau p^{\sigma}})_{p=0}^\infty$, which are not strongly regular, but for which the sequence $\hM$ is regular as long as $1<\sigma<2$. These sequences have been considered in a series of papers by S. Pilipovi{\'c}, N. Teofanov and F. Tomi{\'c}~\cite{ptt15,ptt16,ptt,ptt21}, inducing ultradifferentiable spaces of so-called extended Gevrey regularity. However, it is important to note that the estimates they impose on the elements of those classes are different, as their derivatives are controlled by expressions of the form $CA^{p^\sigma}p^{\tau p^{\sigma}}$ instead. This change of the geometric factor $A^{p}$ into the faster $A^{p^\sigma}$ entails some enhancements in the stability of the classes, but it is not well suited for our purposes, so we attach to the standard framework in ultraholomorphic classes. Our technique is similar to the one employed by S. Malek and the second and third authors~\cite{lastramaleksanzJMAA12} for reproving Thilliez's results, and by the third author in~\cite[Theorem\ 6.1]{SanzFlatProxOrder} for closing the surjectivity problem for strongly regular sequences admitting a nonzero proximate order. It rests on the construction of suitable kernel functions allowing us to define formal Borel- and truncated Laplace-like transforms, in terms of which the solution is provided. These kernels are constructed using the Lambert function and it is possible to prove that they have good behavior thanks to its properties. Regrettably, the multivalued character of the Lambert function, that appears in our procedure in a crucial way, causes some difficulties in order to reason in sectors of opening greater than $2\pi$; for such openings, although surjectivity of the Borel map is known to hold, our constructive proof is not available.

\section{Preliminaries}

\subsection{Notation} Throughout the paper, we use the following notation:  $\N:=\{1,2,\dots \}$, $\N_0=\N\cup\{0\}$, $\C[[z]]$ is the space of formal power series in $z$ with complex coefficients.
We set $\mathcal{R}$ for the Riemann surface of the logarithm, and the notation $z=|z|e^{i\theta}$ refers to the element $(|z|,\theta)\in(0,\infty)\times\R$. For $\gamma>0$ we consider the unbounded sectors
$$S_{\gamma}:=\{z\in \mathcal{R} : |\arg(z)|<\frac{\gamma\pi}{2}\}.$$
For every $0<\gamma<2$ we identify $S_\gamma$ with the corresponding sector in $\C$. For $R>0$, $\gamma>0$, we write
$$L_{R,\gamma}:=\{z\in \mathcal{R} :|z|>R,|\text{arg} (z)|<\frac{\gamma\pi}{2}\}.$$

Given two unbounded sectors $T,S$, we write $T\prec S$ in the case that $\overline{T}\subseteq S\cup\{0\}$. For every $R>0$ and $z_0\in\C$, $D(z_0,R)$ stands for the open disc in $\C$ centered at $z_0$ and radius$~R$.

$\Log(\cdot)$ stands for the principal branch of the complex logarithm, i.e., $\hbox{arg}(\Log(z))\in(-\pi,\pi)$ for every $z\in\C\setminus(-\infty,0]$. We write $\ln(\cdot)$ for the natural logarithm, and for $x>0$, $\ln^+(x)$ denotes the value $\max(\ln(x),0)$.

\subsection{Weight sequences}\label{subsec.WeightSeq} In what follows, $\M=(M_p)_{p\in\N_0}$ will stand for a sequence of positive real numbers with $M_0=1$. The following properties will play a role in this paper:

\begin{itemize}
\item[(i)]  $\M$ is \emph{logarithmically convex} (for short, (lc)) if
$$M_{p}^{2}\le M_{p-1}M_{p+1},\qquad p\in\N.$$%
\item[(ii)] $\M$ is \emph{stable under differential operators} or satisfies the \emph{derivation closedness condition} (briefly, (dc)) if there exists $D>0$ such that $$M_{p+1}\leq D^{p+1} M_{p},  \qquad p\in\N_{0}. $$
\item[(iii)]  $\M$ is of, or has, \emph{moderate growth} (briefly, (mg)) whenever there exists $A>0$ such that
$$M_{p+q}\le A^{p+q}M_{p}M_{q},\qquad p,q\in\N_0.$$
\item[(iv)]  $\M$ satisfies the condition (snq) if there exists $B>0$ such that
$$
\sum^\infty_{q= p}\frac{M_{q}}{(q+1)M_{q+1}}\le B\frac{M_{p}}{M_{p+1}},\qquad p\in\N_0.$$
\end{itemize}

In the classical work of H.~Komatsu~\cite{komatsu}, the properties (lc), (dc) and (mg) are denoted by $(M.1)$, $(M.2)'$ and $(M.2)$, respectively. The property (snq) for $\M$ is precisely the property $(M.3)$ when imposed on the sequence $\widehat{\M}:=(p!M_p)_{p=0}^\infty$.

Obviously, (mg) implies (dc).

If $\M$ is (lc), it is well-known that $((M_p)^{1/p})_{p\in\N}$ is nondecreasing; if moreover $(M_p)^{1/p}\uparrow\infty$, we say $\M$ is a \emph{weight sequence}.
Following E.~M.~Dyn'kin~\cite{Dynkin80}, if $\M$ is a weight sequence and satisfies (dc), we say $\hM$ is \emph{regular}. According to V.~Thilliez~\cite{Thilliez03}, if $\M$ satisfies (lc), (mg) and (snq), we say $\M$ is \emph{strongly regular}; in this case $\M$ is a weight sequence, and the corresponding $\hM$ is regular.

The index $\gamma(\M)$ was introduced by V.~Thilliez~\cite[Sect.\ 1.3]{Thilliez03} for strongly regular sequences $\M$, in which case it is a positive real number, but his definition makes sense for (lc) sequences.

A sequence $(c_p)_{p\in\N_0}$ is \emph{almost increasing} if there exists $a>0$ such that for every $p\in\N_0$ we have that $c_p\leq a c_q $ for every $ q\geq p$; then, it was proved in~\cite{JimenezSanzSRSPO,JimenezSanzSchindlIndex} that for any weight sequence $\M$ one has
\begin{equation}\label{equa.indice.gammaM.casicrec}
\gamma(\bM)=\sup\{\gamma>0: \big(\frac{M_{p+1}}{M_p(p+1)^\gamma}\big)_{p\in\N_0}\hbox{ is almost increasing} \}\in[0,\infty].
\end{equation}

Whenever $\widehat{\M}$ is (lc) we have (see~\cite[Ch.~2]{PhDJimenez} and \cite[Cor.~3.13]{JimenezSanzSchindlIndex}) that
$\gamma(\M)>0$ if and only if $\M$ is (snq).

Two sequences $\M=(M_{p})_{p\in\N_0}$ and $\L=(L_{p})_{p\in\N_0}$ of positive real numbers are said to be \emph{equivalent}, and we write $\M\approx\L$, if there exist positive constants $A,B$ such that
$$A^pM_p\le L_p\le B^pM_p,\qquad p\in\N_0.$$

\subsection{Ultraholomorphic classes}\label{subsec.UltraHolClasses}

In this section $S$ is a sector and $\M$ a sequence.
We say a holomorphic function   $f:S\to\C$ admits the formal power series $\widehat{f}=\sum_{p=0}^{\infty} c_{p}z^{p}\in\C[[z]]$  as its \emph{uniform $\{\M\}$-asymptotic expansion in $S$ (of type $1/A$ for some $A>0$)} if there exists $C>0$ such that for every $p\in\N_0$, one has
\begin{equation}\Big|f(z)-\sum_{p=0}^{n-1}c_pz^p \Big|\le CA^nM_{n}|z|^n,\qquad z\in S.\label{desarasintunifo}
\end{equation}
In this case we write $f\sim_{\{\M\},A}^u\widehat{f}$ in $S$, and $\widetilde{\mathcal{A}}^u_{\{\M\},A}(S)$ denotes the space of functions admitting uniform $\{\M\}$-asymptotic expansion of type $1/A$ in $S$, endowed with the norm
$$
\left\|f\right\|_{\M,A,\overset{\sim}{u}}:=\sup_{z\in S,n\in\N_{0}}\frac{|f(z)-\sum_{p=0}^{n-1}c_pz^p|}{A^{n}M_{n}|z|^n},
$$
which makes it a Banach space. We set $\widetilde{\mathcal{A}}^u_{\{\M\}}(S)$ for the $(LB)$ space of functions admitting uniform $\{\M\}$-asymptotic expansion in $S$, obtained as the union of the previous classes when $A$ runs over $(0,\infty)$.
When the type needs not be specified, we simply write $f\sim_{\{\M\}}^u\widehat{f}$ in $S$.
Note that, taking $n=0$ in~\eqref{desarasintunifo}, we deduce that every function in $\widetilde{\mathcal{A}}^u_{\{\M\}}(S)$ is a bounded function.

One may accordingly define classes of formal power series
$$\C[[z]]_{\{\M\},A}=\Big\{\widehat{f}=\sum_{p=0}^\infty c_pz^p\in\C[[z]]:\, \left|\,\widehat{f} \,\right|_{\M,A}:=\sup_{p\in\N_{0}}\displaystyle \frac{|c_{p}|}{A^{p}M_{p}}<\infty\Big\}.$$
$(\C[[z]]_{\{\M\},A},\left| \  \right|_{\M,A})$ is a Banach space and we put $\C[[z]]_{\{\M\}}:=\cup_{A>0}\C[[z]]_{\{\M\},A}$, again an $(LB)$ space.

Given $f\in\widetilde{\mathcal{A}}^u_{\{\M\}}(S)$ with $f\sim_{\{\M\}}^u\widehat{f}$, Cauchy's integral formula for the derivatives implies that for every $T\prec S$ there exists $A_T>0$ such that
\begin{equation}\label{eq.BoundsDerivUltraholClass}
\sup_{z\in T,p\in\N_{0}}\frac{|f^{(p)}(z)|}{A_T^{p}p!M_{p}}<\infty,
\end{equation}
and then, by Taylor's formula, for every such $T$ and every $p\in\N_0$ one has
\begin{equation*}
c_p=\lim_{ \genfrac{}{}{0pt}{}{z\to0}{z\in T}} \frac{f^{(p)}(z)}{p!},
\end{equation*}
So, we can set ${f^{(p)}(0)}:=p!c_p$, and it is straightforward that $\widehat{f}\in\C[[z]]_{\{\M\}}$, what makes natural to consider the \textit{asymptotic Borel map}
$$
\widetilde{\mathcal{B}}:\widetilde{\mathcal{A}}^u_{\{\M\}}(S)\longrightarrow \C[[z]]_{\{\M\}}$$
sending a function $f\in\widetilde{\mathcal{A}}^u_{\{\M\}}(S)$ into its uniform $\{\M\}$-asymptotic expansion $\widehat{f}$ in $S$. 
$\widetilde{\mathcal{B}}$ is continuous when considered between these  $(LB)$ spaces, and also when restricted to the Banach spaces with fixed type.

Finally, note that if $\M\approx\L$, then $\widetilde{\mathcal{A}}^u_{\{\M\}}(S)=\widetilde{\mathcal{A}}^u_{\{\L\}}(S)$ and $\C[[z]]_{\{\M\}}=\C[[z]]_{\{\L\}}$, so the corresponding Borel maps are in all cases identical.

\begin{rem}
The inequalities mentioned in~\eqref{eq.BoundsDerivUltraholClass} suggest considering ultraholomorphic classes defined by imposing such estimates on the derivatives of its elements $f$, and the corresponding Borel map sending $f$ into $(f^{(p)}(0))_{p\in\N_0}$, where $f^{(p)}(0)$ is defined as above. Another possibility consists in regarding non-uniform asymptotics, with bounds $C_TA_T^pM_p|z|^p$ for the remainders on every bounded and proper subsector $T$ of $S$. Our forthcoming results can be stated in these frameworks, but we will not enter into details. The interested reader may consult~\cite{JimenezSanzSchindlInjectSurject,JimenezSanzSchindlSurjectDc}.
\end{rem}

\subsection{A family of sequences of rapid growth}

Let $\sigma>1$ and $\tau>0$ be real numbers. We consider the sequence of positive real numbers $\mathbb{M}^{\tau,\sigma}=(M_p^{\tau,\sigma})_{p\in\N_0}$ with $M_p^{\tau,\sigma}:=p^{\tau p^{\sigma}}$ for $p\ge1$ and $M^{\tau,\sigma}_0=1$.

First we recall the essential properties of the sequences $\mathbb{M}^{\tau,\sigma}$.
We refer to \cite{ptt15} for the proof of the next Lemma.

\begin{lemma} Let  $\sigma>1$ and $\tau>0$ be real numbers. Then
the following properties hold:
 \begin{enumerate}
  \item[(i)] $\mathbb{M}^{\tau,\sigma}$ is (lc).
  \item[(ii)] $\overline{(M.2)'}$ For all $q\in\N_{0}$ there exists $C_q\geq 1$ such that
  $$M^{\tau,\sigma}_{p+q}\leq C_q^{p^\sigma} M^{\tau,\sigma}_{p},  \qquad p,q\in\N_{0}. $$
\item[(iii)]  $\overline{(M.2)}$ There exists $A>0$ such that
$$M^{\tau,\sigma}_{p+q}\le A^{p^\sigma +q^\sigma}M^{\tau 2^{\sigma-1},\sigma}_{p}M^{\tau 2^{\sigma-1},\sigma}_{q},\qquad p,q\in\N_0.$$
 \end{enumerate}
\end{lemma}

We note that the notation $\overline{(M.2)'}$, resp. $\overline{(M.2)}$, comes from~\cite{ptt}, while in~\cite{ptt15,ptt16,ptt21} $\widetilde{(M.2)'}$, resp. $\widetilde{(M.2)}$, is used.

These conditions $\overline{(M.2)'}$ and $\overline{(M.2)}$ are different from the classical ones, namely (dc) and (mg), appearing in the literature when dealing with Carleman-like classes (see~\cite{komatsu}). They play a prominent role in the study of the corresponding ultradifferentiable and ultradistributional classes carried out by S. Pilipovi{\'c}, N. Teofanov and F. Tomi{\'c} in~\cite{ptt15,ptt16,ptt,ptt21}, as they allow for a precise control of the flexibility obtained by introducing the two parameter dependence.

Accordingly, they introduce a convenient modification of the classical associated function by setting (see~\cite[Definition~2.1]{ptt})
$$T_{\tau,\sigma,h}(t):= \sup_{p\in\mathbb{N}} \ln^+\left(\frac{h^{p^{\sigma}}t^{p}}{M_{p}^{\tau,\sigma}}\right),\qquad t>0,$$
and $T_{\tau,\sigma,h}(0)=0$. We write $T_{\tau,\sigma}$ for the function $T_{\tau,\sigma,1}$, which is the standard function associated with the sequence $\M^{\tau,\sigma}$ in the literature (see~\cite[Definition~3.1]{komatsu}).

\begin{rem}\label{rem.RapidGrowth}

We note that, from the expression~\eqref{equa.indice.gammaM.casicrec}, it is easy to check that $\gamma(\M^{\tau,\sigma})=+\infty$ for every $\sigma>1$, so we say these sequences are of rapid growth. The terminology comes from the theory of regular variation, see~\cite{bgt} for the case of functions and~\cite{DjurcicKocinacZizovic07,JimenezSanzSchindlIndex} for some extensions of the theory for the case of sequences. Although there are different definitions of the notion of rapid variation, they all coincide when we restrict ourselves to nondecreasing sequences, see \cite{DjurcicElezKocinac15}.  Then one can say that a nondecreasing  sequence is of rapid variation if its lower Matuszewska index is $+\infty$. It turns out~\cite[Section~3.3]{JimenezSanzSchindlIndex} that, for a general sequence of positive real numbers $\M=(M_p)_{p\in\N_0}$, the index $\gamma(\M)$ equals the lower Matuszewska index of the sequence of quotients $(M_{p+1}/M_p)_{p\in\N_0}$, which is nondecreasing if $\M$ is (lc); so, the sequence $(M^{\tau,\sigma}_{p+1}/M^{\tau,\sigma}_p)_{p\in\N_0}$ is of rapid variation, what justifies the previous nomenclature.

As mentioned before, the rapid growth of $\M^{\tau,\sigma}$ implies (snq) is satisfied, while (mg) is not, hence $\M^{\tau,\sigma}$ is not strongly regular for any $\sigma>1$. However, whenever $\sigma$ belongs to the interval $(1,2)$, the sequence $(M_p^{\tau,\sigma})_{p\in\N_0}$ not only satisfies $\overline{(M.2)'}$  but also  the property (dc), i.e., there exists $C>0$ such that $M_{p+1}^{\tau,\sigma}\le C^{p+1}M_{p}^{\tau,\sigma}$. Indeed, for logarithmically convex sequences $(M_p)_{p\in\N_0}$, (dc) is equivalent to the condition $\log(M_p)=O(p^2)$, $p\to\infty$, which clearly holds for $\M^{\tau,\sigma}$ if, and only if, $\sigma\in(1,2)$. Hence, for such values $\widehat{\M}^{\tau,\sigma}$ is regular. This fact will make a difference with respect to the general case. Moreover, it implies the following well-known estimate: if one defines the auxiliary function
$$
h_{\tau,\sigma}(t):=e^{-T_{\tau,\sigma}(1/t)}=\inf_{p\ge 0}M_{p}^{\tau,\sigma}t^p,\quad t>0;\ h_{\tau,\sigma}(0)=0,
$$
for every $t>0$ one has
\begin{equation}
h_{\tau,\sigma}(t)\le
\inf_{p\ge 1}M_{p}^{\tau,\sigma}t^p=\inf_{p\ge 0}M_{p+1}^{\tau,\sigma}t^{p+1}\le Ct\inf_{p\ge 0}M_{p}^{\tau,\sigma}(Ct)^p=Ct h_{\tau,\sigma}(Ct).\label{eq.Estimate_h_tausigma_dc}
\end{equation}
\end{rem}

\section{On the Lambert function $W$}

In this section, and for the sake of completeness, we give some information about the Lambert function $W$ and its main properties to be considered in the present study. We refer to~\cite{corless} for further details.

The Lambert function $W$ is defined as the complex function which satisfies
\begin{equation}\label{eW}
W(z)e^{W(z)}=z.
\end{equation}
It holds that the Lambert function $W$ is a multivalued function which splits the $w=W(z)$ plane into an infinite number of regions. The principal branch of the Lambert function $W$, usually denoted by $W_0$, is defined in $\C\setminus (-\infty,-e^{-1}]$, and the curve in $\R^2$ defining the boundary of its image set $W_0(\C\setminus (-\infty,-e^{-1}])$ is given by
$$\{(-t\cot(t),t)\in\R^2:-\pi<t<\pi\}.$$
$W_0$ is a holomorphic and bijective map when restricted to the previous domains, with $W_0(0)=0$, being the origin the only value of its domain in which $W$ vanishes. Hereinafter, we will restrict our reasonings to the function $W_0$, so we will write $W$ instead of $W_0$ from now on.

\begin{figure}
	\centering
		\includegraphics[width=0.5\textwidth]{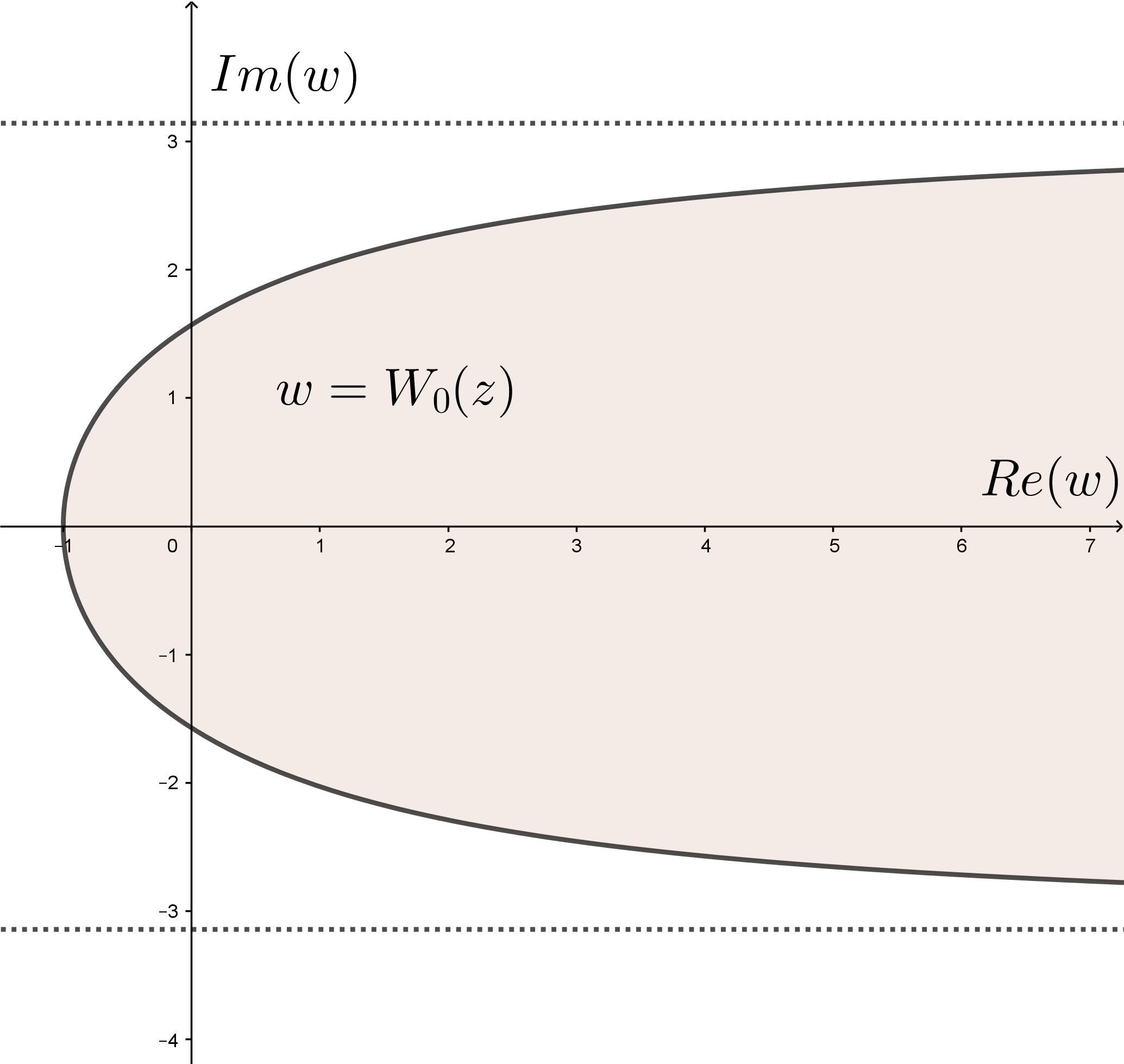}
		\caption{Principal branch of the Lambert function $W$}
\end{figure}

\begin{prop}\label{propW}
The following properties hold for the principal branch of the Lambert function $W$:
\begin{enumerate}
\item[(a)] $W(0)=0$. More precisely, $z=0$ is a zero of order 1, and indeed $\displaystyle\lim_{z\to 0}W(z)/z=1$.
\item[(b)] $W(x)\in\R$ for $x>-e^{-1}$, and $\lim_{x\to\infty}W(x)=+\infty$. In addition to this, one has that $W(x)>0$ for $x>0$.
\item[(c)] For every $z\in\C\setminus(-\infty,-e^{-1}]$ one has
$$W'(z)=e^{-W(z)}\frac{1}{1+W(z)}.$$
\item[(d)] For every $z\in\C\setminus(-\infty,-e^{-1}]$ it holds that
$$\Re(z)=e^{\beta_1(z)}(\beta_1(z)\cos(\beta_2(z))-\beta_2(z)\sin(\beta_2(z))),$$
$$\Im(z)=e^{\beta_1(z)}(\beta_2(z)\cos(\beta_2(z))+\beta_1(z)\sin(\beta_2(z))),$$
where $\beta_1(z)=\Re(W(z))$ and $\beta_2(z)=\Im(W(z))$.
\item[(e)] Let $R>0$ and $\alpha\in (0,2)$. For $L_{R,\alpha}=S_{\alpha}\cap (\C\setminus D(0,R))$ one has
\begin{equation}\label{eq.imageunderW_kofasector_bis}
W(L_{R,\alpha})\subseteq\{\xi+i\eta \in\mathbb{C}: \exp(\xi)\sqrt{\xi^2+\eta^2}\ge R\}.
\end{equation}
\end{enumerate}
\end{prop}
\begin{proof1}
(a) is a direct consequence of (\ref{eW}), from where the Taylor expansion of the Lambert function $W$ at the origin of the complex plane can be obtained, see~\cite{corless}. (b) also stems from (\ref{eW}), taking into account that $W$, when restricted to the interval $(-1,\infty)$, is the inverse function of $x\mapsto xe^x$. Differentiating (\ref{eW}) one arrives at an expression equivalent to (c). Observe that $W(z)\neq-1$ for $z\in\C\setminus(-\infty,-e^{-1}]$. (d) and (e) are obtained by splitting (\ref{eW}) into its real and imaginary parts.
\end{proof1}

\begin{rem}\label{remWslow}
It is worth remarking that property (e) in the previous proposition guarantees that for every $R>0$ and $\alpha\in(0,2)$,
\begin{equation}\label{e119}
\lim_{|z|\to\infty}\frac{zW'(z)}{W(z)}=\lim_{|z|\to\infty}\frac{1}{1+W(z)}=0,
\end{equation}
uniformly in $L_{R,\alpha}$. Indeed, given $M>0$ we may choose $R_0\ge Me^M$. For every $z\in L_{R_0,\alpha}$
put $W(z)=\xi +i\eta$. There are two possibilities: if $\xi\ge M$ then obviously $|W(z)|\ge M$; on the contrary, \eqref{eq.imageunderW_kofasector_bis} gives $|W(z)|=(\xi^2+\eta^2)^{1/2}\ge R_0e^{-\xi}>Me^Me^{-M}=M$, as desired.

Property (\ref{e119}) states that $W$ is slowly varying (see~\cite[Theorem~$A.1.2.(a)$]{bgt}) in $L_{R,\alpha}$. This statement is also true for any other choice of the branch of the Lambert function.
\end{rem}

\section{Kernel functions for sequences of rapid growth}

In this section, we construct a kernel function associated with the sequence $\M^{\tau,\sigma}$ which will be useful for our purposes.

\begin{defi}\label{defib}
For $\sigma>1$ and $\tau>0$, we define 
\begin{equation}\label{e177}
a_{\tau,\sigma}:= \left(\frac{\sigma-1}{\tau\sigma}\right)^{\frac{1}{\sigma-1}}, \qquad b_{\tau,\sigma}:=e^{\frac{\sigma-1}{\sigma}} \frac{\sigma-1}{\tau\sigma}.
\end{equation}
and 
$$
e_{\tau,\sigma}(z):=z\exp\left( -a_{\tau,\sigma} W^{-\frac{1}{\sigma-1}}(b_{\tau,\sigma}\Log(z+1)) \Log^{\frac{\sigma}{\sigma-1}}(z+1)\right)
$$
for every $z\in\mathbb{C}\setminus(-\infty,0]$. 
\end{defi}

For the sake of brevity, we write
\begin{equation}\label{e351}
 g(z):=W^{-\frac{1}{\sigma-1}}(b_{\tau,\sigma}\Log(z)) \Log^{\frac{\sigma}{\sigma-1}}(z).
\end{equation}
so $e_{\tau,\sigma}(z)=z\exp(-a_{\tau,\sigma}g(z+1))$.

The powers of the logarithm and of the Lambert function appearing in $e_{\tau,\sigma}$ and $g$ have their principal values, i.e., they are the ones associated with the principal branch of the logarithm.

\begin{prop}\label{prop347}
$e_{\tau,\sigma}(z)$ is a holomorphic function in $\mathbb{C}\setminus(-\infty,0]$.
\end{prop}
\begin{proof1}
On the one hand, if $z\in\mathbb{C}\setminus(-\infty,0]$, then $z+1\in\mathbb{C}\setminus(-\infty,1]$, and so
$\Log(z+1)=\ln(|z+1|)+i\arg(z+1)\notin(-\infty,0]$. Since $b_{\tau,\sigma}\Log(z+1)$ differs from $\Log(z+1)$ by a positive factor, this fact ensures that the maps
\begin{align*}
z&\mapsto \Log^{\frac{\sigma}{\sigma-1}}(z+1)= \exp\Big(\frac{\sigma}{\sigma-1}\Log\big(\Log(z+1)\big)\Big),\\
z&\mapsto W(b_{\tau,\sigma}\Log(z+1)),
\end{align*}
are both holomorphic in $\mathbb{C}\setminus(-\infty,0]$. Finally, observe that the principal branch of the Lambert function $W(w)$ only takes values in $(-\infty,0]$ when $w$ runs over the real interval $(-1/e,0]$. Since $w=b_{\tau,\sigma}\Log(z+1)\notin(-\infty,0]$ whenever $z\in\mathbb{C}\setminus(-\infty,0]$, the composition of $z\mapsto W(b_{\tau,\sigma}\Log(z+1))$ with the principal branch of the logarithm is again holomorphic, and so is
$$
z\mapsto W^{-\frac{1}{\sigma-1}}(b_{\tau,\sigma}\Log(z+1)),
$$
what leads to the conclusion.
\end{proof1}

\begin{rem}
There is some freedom in the choice of $e_{\tau,\sigma}$. The factor $z$ may be changed into any $z^{\alpha}$ for some positive real number $\a$ (so that the assertion $(i)$ in Lemma~\ref{lema.Properties.e_h} holds true), where the principal branch of the power is considered. Our choice tries to make the following computations simpler.
\end{rem}

\begin{rem}\label{rem.e_h_boundedat0}
Taking into account property (a) of Proposition~\ref{propW}, one has that
\begin{align*}
\lim_{z\to 0}\frac{e_{\tau,\sigma}(z)}{z}&=\lim_{z\to 0}\exp\left(-\left(\frac{\sigma-1}{\tau\sigma}\right)^{\frac{1}{\sigma-1}} W^{-\frac{1}{\sigma-1}}(b_{\tau,\sigma}\Log(z+1))\Log^{\frac{\sigma}{\sigma-1}}(z+1)\right)\\
&=\lim_{z\to 0}\exp\left(-\left(\frac{\sigma-1}{\tau\sigma}\right)^{\frac{1}{\sigma-1}} \frac{\Log^{\frac{\sigma}{\sigma-1}}(z+1)}{(b_{\tau,\sigma}\Log(z+1))^{\frac{1}{\sigma-1}}}\right)\\
&=\lim_{z\to 0}\exp\left(-\left(\frac{\sigma-1}{\tau\sigma}\right)^{\frac{1}{\sigma-1}} \frac{\Log^{\frac{\sigma}{\sigma-1}}(z+1)}{e^{\frac{1}{\sigma}} \left(\frac{\sigma-1}{\tau\sigma}\right)^{\frac{1}{\sigma-1}} \Log^{\frac{1}{\sigma-1}}(z+1)}\right)\\
&=\lim_{z\to 0}\exp\left(-e^{-\frac{1}{\sigma}}
\Log(z+1)\right)=1.
\end{align*}
Since the function $e_{\tau,\sigma}(z)/z$ never vanishes, we deduce that it will remain bounded, and bounded away from 0, in bounded proper sectors of $\mathbb{C}\setminus(-\infty,0]$ with vertex at 0.
\end{rem}

\begin{lemma}\label{lema.g_hincreasing}
The restriction of the function $g$ in \eqref{e351} to $(1,\infty)$, given by
\begin{equation}\label{eq.defi.g_h_reals}
g(x)=W^{-\frac{1}{\sigma-1}}(b_{\tau,\sigma}\ln(x)) \ln^{\frac{\sigma}{\sigma-1}}(x),
\end{equation}
is positive and strictly increasing. Moreover, $\lim_{x\to\infty}g'(x)=0$.
\end{lemma}

\begin{proof1}
We clearly have $b_{\tau,\sigma}\ln(x)>0$ for every $x>1$, and $W$ is positive for positive values of its argument (see (b) in Proposition~\ref{propW}), hence $g(x)>0$ for every $x>1$. Property (c) of Proposition~\ref{propW} yields
\begin{equation}\label{eq.derivW}
W'(x)=\frac{1}{x+e^{W(x)}}=\frac{W(x)}{x(W(x)+1)}.
\end{equation}
So, applying~\eqref{eq.derivW} we can write
\begin{align*}
g'(x)&=-\frac{1}{\sigma-1}W^{-\frac{\sigma}{\sigma-1}}(b_{\tau,\sigma}\ln(x)) W'(b_{\tau,\sigma}\ln(x))b_{\tau,\sigma}\frac{1}{x}\ln^{\frac{\sigma}{\sigma-1}}(x)\\
&+W^{-\frac{1}{\sigma-1}}(b_{\tau,\sigma}\ln(x))\frac{\sigma}{\sigma-1} \ln^{\frac{1}{\sigma-1}}(x)\frac{1}{x}\\
&=-\frac{1}{\sigma-1}W^{-\frac{1}{\sigma-1}}(b_{\tau,\sigma}\ln(x)) \frac{1}{W(b_{\tau,\sigma}\ln(x))+1}\frac{1}{x}\ln^{\frac{1}{\sigma-1}}(x)\\
&+W^{-\frac{1}{\sigma-1}}(b_{\tau,\sigma}\ln(x))\frac{\sigma}{\sigma-1} \ln^{\frac{1}{\sigma-1}}(x)\frac{1}{x}\\
&=\frac{1}{\sigma-1}W^{-\frac{1}{\sigma-1}}(b_{\tau,\sigma}\ln(x)) \ln^{\frac{1}{\sigma-1}}(x)\frac{1}{x} \left(\sigma-\frac{1}{W(b_{\tau,\sigma}\ln(x))+1}\right)>0,
\end{align*}
since $x>1$ and $\sigma>1$. So, $g$ is strictly increasing in $(1,\infty)$. If we observe now that $\lim_{x\to\infty}b_{\tau,\sigma}\ln(x)=+\infty$ and (b) in Proposition~\ref{propW}, it is clear from the last expression that $\lim_{x\to\infty}g'(x)=0$.
\end{proof1}

The next technical result provides information on the function $e_{\tau,\sigma}$ with respect to the standard properties satisfied by kernel functions for summability, as described in~\cite{lastramaleksanzJMAA15}, and which were inspired by the work of W. Balser~\cite[Section~5.5]{balser} in the context of Gevrey classes.

\begin{lemma}\label{lema.Properties.e_h}
The function $e_{\tau,\sigma}(\cdot)$ enjoys the following properties:
\begin{itemize}
\item[1.] $z^{-1}e_{\tau,\sigma}(z)$ is integrable at the origin, it is to say, for every $t_0>0$ and $\theta\in(-\pi,\pi)$, the integral
$$\int_0^{t_0}t^{-1}|e_{\tau,\sigma}(te^{i\theta})|dt<\infty.$$
\item[2.] For every $x>0$ the value of $e_{\tau,\sigma}(x)$ is positive real.
\item[3.] For every $T\prec S_{2}$ there exist $C_1,K_1,C_2,K_2>0$ (depending on $\tau,\sigma$) such that
\begin{equation}\label{eq.Bounds_e_h_subsectors_real_axis}
C_1e_{\tau,\sigma}(K_1|z|)\le |e_{\tau,\sigma}(z)|\le C_2e_{\tau,\sigma}(K_2|z|),\quad z\in T.
\end{equation}
\end{itemize}
\end{lemma}

\begin{proof1}
\noindent 1. By the very definition of the function $e_{\tau,\sigma}$, it is evident that the function $z^{-1}e_{\tau,\sigma}(z)$ is measurable and bounded in every segment towards the origin, and so it is also integrable.

\noindent 2. Recall that the function $g$ in (\ref{eq.defi.g_h_reals})  is the restriction of the function defined in (\ref{e351}) to $(1,\infty)$. Then it suffices to observe that the value $g(x+1)$ is real for every $x>0$, see Lemma~\ref{lema.g_hincreasing}.

\noindent 3. Throughout this part of the proof, we make use of property (e) in Proposition~\ref{propW} to take control of the image of the elements contained in a sector of the complex plane which stay away from the origin via Lambert function.

 The constants $a_{\tau,\sigma}$ and $b_{\tau,\sigma}$, defined in (\ref{e177}), will be named $a$ and $b$ throughout this proof for brevity.
For $\omega\neq 0$ one can apply (\ref{eW}) to $z=b\omega$ from where we deduce
$$
\frac{\omega^{\frac{\sigma}{\sigma-1}}}{W^{\frac{1}{\sigma-1}}(b\omega)}=
\frac{1}{b^{1/(\sigma-1)}}\omega e^{\frac{1}{\sigma-1}W(b\omega)}.
$$
So, putting $\omega=\Log(z+1)$ we obtain
$$
g(z+1)=\frac{1}{b^{1/(\sigma-1)}}\Log(z+1) e^{\frac{1}{\sigma-1}W(b\Log(z+1))}.
$$
According to this expression, we will start by showing that given $\beta\in(0,2)$, there exist $R,K>0$ such that for every $z\in S_{\beta}$ with $|z|\ge R$ one has
$$
\Re\Big(\Log(z+1) e^{\frac{1}{\sigma-1}W(b\Log(z+1))}\Big)\le
\ln(K|z|+1) e^{\frac{1}{\sigma-1}W(b\ln(K|z|+1))};
$$
equivalently, if we write
\begin{equation}\label{eq.RealImaginaryPartsLambert}
W(b\Log(z+1))=\xi(z)+i\eta(z),\quad \xi(z),\eta(z)\in\R,
\end{equation}
one has
\begin{multline}\label{eq.IneqFineAdjustComplexLambert}
\ln(|z+1|) e^{\frac{1}{\sigma-1}\xi(z)} \cos\big(\frac{1}{\sigma-1}\eta(z)\big)-
\arg(z+1)e^{\frac{1}{\sigma-1}\xi(z)} \sin\big(\frac{1}{\sigma-1}\eta(z)\big)\\
\le
\ln(K|z|+1) e^{\frac{1}{\sigma-1}W(b\ln(K|z|+1))},
\end{multline}
where $W$ is the principal branch for the Lambert function. By taking exponentials in the last inequality, we obtain the first inequality in~\eqref{eq.Bounds_e_h_subsectors_real_axis}, if we take into account that the function
$$
\frac{|e_{\tau,\sigma}(z)|}{e_{\tau,\sigma}(K|z|)}= \frac{\exp(-a\Re(g(z+1)))}{K\exp(-ag(K|z|+1))}
$$
is bounded above, and bounded away from zero, in bounded proper subsectors of $\C\setminus(-\infty,0]$ (see Remark~\ref{rem.e_h_boundedat0}).

Observe that if $R>2$ and $z\in S_{\beta}$ with $|z|\ge R$, it is clear that $z+1\in S_{\beta}$ and $|z+1|\ge R-1>1$, so that we clearly have
\begin{multline}\label{eq.FirstIneqWlogzWlogKz}
\ln(|z+1|) e^{\frac{1}{\sigma-1}\xi(z)} \cos\big(\frac{1}{\sigma-1}\eta(z)\big)-
\arg(z+1)e^{\frac{1}{\sigma-1}\xi(z)} \sin\big(\frac{1}{\sigma-1}\eta(z)\big)\\
\le e^{\frac{1}{\sigma-1}\xi(z)}\Big(\ln(|z+1|)+\frac{\pi\beta}{2}\Big).
\end{multline}
Since $\lim_{x\to +\infty}(\ln(Kx+1)-\ln(x+1))=\ln(K)$, if we choose $K>e^{\pi\beta/2}$ there exists $R>2$ such that for every $z\in S_{\beta}$ with $|z|\ge R$ one has
\begin{equation}\label{eq.IneqLogarithmModulez}
\ln(|z+1|)+\frac{\pi\beta}{2}\le\ln(|z|+1)+\frac{\pi\beta}{2}\le\ln(K|z|+1).
\end{equation}
From this fact we also deduce that for such $z$,
\begin{equation*}
\ln^2(K|z|+1)-\ln^2(|z|+1)= \ln\big((K|z|+1)(|z|+1)\big)\ln\Big(\frac{K|z|+1}{|z|+1}\Big)
\ge\Big(\frac{\pi\beta}{2}\Big)^2.
\end{equation*}
Hence, if we write $W(b\ln(K|z|+1))=\xi_1(z)\in\R$, and we take into account~\eqref{eq.RealImaginaryPartsLambert} and~\eqref{eW}, we can obtain that
\begin{align*}
\xi(z)e^{\xi(z)}&\le e^{\xi(z)}(\xi^2(z)+\eta^2(z))^{1/2}=|b\Log(z+1)|\\
&=b\big(\ln^2(|z+1|)+\arg^2(z+1)\big)^{1/2}\le b\Big(\ln^2(|z|+1)+\Big(\frac{\pi\beta}{2}\Big)^2\Big)^{1/2}\\
&\le b\ln(K|z|+1)=\xi_1(z)e^{\xi_1(z)}.
\end{align*}
Since the map $\xi\in(-1,\infty)\mapsto\xi e^{\xi}$ is strictly increasing, we have
\begin{equation}\label{eq.IneqRealPartLambertandLambert}
\xi(z)\le\xi_1(z)=W(b\ln(K|z|+1)).
\end{equation}
Gathering \eqref{eq.FirstIneqWlogzWlogKz},~\eqref{eq.IneqLogarithmModulez} and~\eqref{eq.IneqRealPartLambertandLambert}, we conclude that~\eqref{eq.IneqFineAdjustComplexLambert} is satisfied, and so we get the first inequality in~\eqref{eq.Bounds_e_h_subsectors_real_axis} for a suitably enlarged constant $C_2\ge 1$.

We turn now to the second inequality in~\eqref{eq.Bounds_e_h_subsectors_real_axis}. Reasoning as before, it suffices to prove that given $\beta\in(0,2)$, there exist $R,H>0$ such that for every $z\in S_{\beta}$ with $|z|\ge R$ one has
\begin{multline}\label{eq.IneqFineAdjustBelowComplexLambert}
e^{\frac{1}{\sigma-1}\xi(z)}\Big(\ln(|z+1|) \cos\big(\frac{1}{\sigma-1}\eta(z)\big)-
\arg(z+1)\sin\big(\frac{1}{\sigma-1}\eta(z)\big)\Big)\\
\ge
\ln\big(\frac{|z|}{H}+1\big) e^{\frac{1}{\sigma-1}W\big(b\ln\big(\frac{|z|}{H}+1\big)\big)}.
\end{multline}
In view of Remark~\ref{remWslow} we derive that $W(b\Log(z+1))$ tends to infinity as $z$ does so in $S_{\beta}$. Since $|\eta(z)|\le\pi$ for the principal branch of the Lambert function, necessarily $\xi(z)=\Re(W(b\Log(z+1)))$ tends to $+\infty$ as $z\to\infty$ in $S_{\beta}$. Now, in view of property (d) in Proposition~\ref{propW} it is straightforward that
\begin{equation}\label{eq.EqualityImaginaryPartsLambert}
\Im(b\Log(z+1))=e^{\xi(z)}\big(\eta(z)\cos(\eta(z))+ \xi(z)\sin(\eta(z))\big).
\end{equation}
As the left-hand side is bounded by $b\beta\pi/2$ for $z$ in $S_{\beta}$, we deduce that
$$
\lim_{z\to\infty,\ z\in S_{\beta}}\big(\eta(z)\cos(\eta(z))+\xi(z)\sin(\eta(z))\big)=0.
$$
Again, $\eta(z)\cos(\eta(z))$ is bounded, so the only possibility is that $\lim_{z\to\infty,\ z\in S_{\beta}}\sin(\eta(z))=0$, and hence $\lim_{z\to\infty,\ z\in S_{\beta}}\eta(z)=0$. Here we have used the fact that, due to a connectedness argument, $W(b\Log(z+1))$ cannot approach the boundary of the image set of the principal branch of the Lambert function, so that $\eta(z)$ cannot tend to either $\pi$ or $-\pi$. Going back to~\eqref{eq.EqualityImaginaryPartsLambert}, we get that, as $z$ tends to $\infty$, $\eta(z)$ is equivalent to
$$\sin(\eta(z))=\frac{b\arg(z+1)e^{-\xi(z)}-\eta(z)\cos(\eta(z))}{\xi(z)}.$$
This means there exists a function $\varepsilon(z)$, tending to 1 as $z\to\infty$, such that
$$
\varepsilon(z)\eta(z)= \frac{b\arg(z+1)e^{-\xi(z)}-\eta(z)\cos(\eta(z))}{\xi(z)},
$$
and so
\begin{align}\label{eq.AsymptoticsEtaz}
\eta(z)&= \frac{b\arg(z+1)e^{-\xi(z)}}{\varepsilon(z)\xi(z)+\cos(\eta(z))}
\sim  \frac{b\arg(z+1)}{\xi(z)e^{\xi(z)}},\quad z\to\infty,\\
e^{\xi(z)}\eta(z)&\sim \frac{b\arg(z+1)}{\xi(z)}, \quad z\to\infty.
\label{eq.AsymptoticsEtazExpEtaz}
\end{align}
Again from the definition of the Lambert function, it is straightforward that
$$b\ln(|z+1|)=\Re(b\Log(z+1))=e^{\xi(z)}\big(\xi(z)\cos(\eta(z))- \eta(z)\sin(\eta(z))\big),$$
from where we obtain that
\begin{equation}\label{eq.AsymptoticsbLnz+1}
b\ln(|z+1|)\sim \xi(z)e^{\xi(z)},\quad z\to\infty,\ z\in S_{\beta}.
\end{equation}
Observe that
\begin{equation}\label{eq.RewriteLnz+1Cos}
\ln(|z+1|) \cos\big(\frac{1}{\sigma-1}\eta(z)\big)=
\ln(|z+1|) -\ln(|z+1|) \frac{\sin^2\big(\frac{1}{\sigma-1}\eta(z)\big)}{1+ \cos\big(\frac{1}{\sigma-1}\eta(z)\big)}
\end{equation}
and then, because of~\eqref{eq.AsymptoticsbLnz+1} and~\eqref{eq.AsymptoticsEtaz},
\begin{equation}\label{eq.AsymptoticsLnz+1sinSquare}
\ln(|z+1|)\sin^2\big(\frac{1}{\sigma-1}\eta(z)\big)\sim \frac{\xi(z)e^{\xi(z)}b^2\arg^2(z+1)}{b(\sigma-1)^2 \big(\xi(z)e^{\xi(z)}\big)^2}\to 0 \ \text{as }z\to\infty.
\end{equation}
Hence, going back to~\eqref{eq.IneqFineAdjustBelowComplexLambert}, on the one hand there exists $R_1>2$ such that for every $z\in S_{\beta}$ with $|z|\ge R_1$ we have, by virtue of~\eqref{eq.RewriteLnz+1Cos} and~\eqref{eq.AsymptoticsLnz+1sinSquare},
$$
\ln(|z+1|) \cos\big(\frac{1}{\sigma-1}\eta(z)\big)-
\arg(z+1)\sin\big(\frac{1}{\sigma-1}\eta(z)\big)\ge
\ln(|z+1|)-\frac{1}{b}\ge \ln(|z|-1)-\frac{1}{b}.
$$
Since, for $H>0$, $\lim_{x\to +\infty}(\ln(x-1)-\ln(x/H+1))=\ln(H)$, if we choose $H>e^{1/b}$ there exists $R_2>R_1$ such that for $x\ge R_2$ one has
\begin{equation}\label{eq.Lnx-1LnxH+1}
\ln(x-1)-\ln\big(\frac{x}{H}+1\big)>\frac{1}{b},
\end{equation}
and so for every $z\in S_{\beta}$ with $|z|\ge R_2$ one gets
\begin{equation}\label{eq.IneqLogarithmModulez_bis}
\ln(|z+1|) \cos\big(\frac{1}{\sigma-1}\eta(z)\big)-
\arg(z+1)\sin\big(\frac{1}{\sigma-1}\eta(z)\big)>
\ln\big(\frac{|z|}{H}+1\big).
\end{equation}
On the other hand, by~\eqref{eq.AsymptoticsEtazExpEtaz} there exists
$R_3>2$ such that $e^{\xi(z)}|\eta(z)|<1$. If we take $z\in S_{\beta}$ with $|z|\ge R_4:=\max(R_3,R_2)$, we deduce that
\begin{align*}
\xi(z)e^{\xi(z)}&=e^{\xi(z)}\Re(W(b\Log(z+1)))\\
&\ge e^{\xi(z)}\big(|W(b\Log(z+1))|-|\Im(W(b\Log(z+1)))|\big)\\
&= b|\Log(z+1)|-e^{\xi(z)}|\eta(z)|\ge b\ln(|z+1|)-1\\
&\ge b\ln(|z|-1)-1\ge b\ln\big(\frac{|z|}{H}+1\big),
\end{align*}
where the last inequality comes from~\eqref{eq.Lnx-1LnxH+1}. Due again to the fact that $\xi\in(-1,\infty)\mapsto\xi e^{\xi}$ is strictly increasing, we have that $\xi(z)\ge W(b\ln(|z|/H+1))$. Together with~\eqref{eq.IneqLogarithmModulez_bis}, this proves~\eqref{eq.IneqFineAdjustBelowComplexLambert}.
\end{proof1}

For every $p\in\mathbb{N}_0$ we define the $p-$th moment associated with the kernel function $e_{\tau,\sigma}$ by
$$m_{\tau,\sigma}(p):=\int_0^\infty t^{p-1}e_{\tau,\sigma}(t)\,dt=
\int_0^\infty t^p\exp\left(-a_{\tau,\sigma} g(t+1)\right)\,dt.$$

In order to compare the sequence of moments and the original sequence $\M^{\tau,\sigma}$, the following result~\cite[Theorem~2.1]{ptt} is important. It is particularized for the case $1<\sigma<2$, which we are considering here; a more general result is available which, however, does not meet our needs, see Remark~\ref{rem.ChangeBoundsSigmaGreater2}.

\begin{lemma}[\cite{ptt}, Theorem 2.1, Remark 2.2]\label{lema.doubleEstimatesSerbian}
Let $1<\sigma<2$.
There exist $A_{\tau,\sigma},\tilde{A}_{\tau,\sigma}>0$ such that
\begin{equation}\label{eq.doubleEstimatesSerbian}
A_{\tau,\sigma}\exp(a_{\tau,\sigma}g(x))\le e^{T_{\tau,\sigma}(x)}\le \tilde{A}_{\tau,\sigma}\exp(a_{\tau,\sigma}g(x)),
\end{equation}
for every $x>1$.
\end{lemma}

\begin{prop}\label{prop.boundsMomentsSigmaEntre1y2}

Suppose $1<\sigma<2$. Consider the sequence of moments $(m_{\tau,\sigma}(p))_{n\ge0}$ associated with the kernel function $e_{\tau,\sigma}(z)$. Then, there exist $B_1,B_2>0$ such that
\begin{equation}\label{eq.boundsMmomentsdc}
B_1^{p}M_p^{\tau,\sigma}\le m_{\tau,\sigma}(p)\le B_2^{p}M_p^{\tau,\sigma},\quad p\in\mathbb{N}_0.
\end{equation}
\end{prop}

\begin{proof1}
On the one hand, for every $s>0$ we may write
\begin{align}
m_{\tau,\sigma}(p)&=\int_0^s t^p\exp\left(-a_{\tau,\sigma}g(t+1)\right)\,dt
+\int_s^\infty \frac{1}{t^2}t^{p+2}\exp\left(-a_{\tau,\sigma}g(t+1)\right)\,dt\nonumber\\
&\le \int_0^s t^p\,dt
+\frac{1}{s}\sup_{t>0}(t+1)^{p+2}\exp\left(-a_{\tau,\sigma}g(t+1)\right)\nonumber\\
&\le \frac{s^{p+1}}{p+1}+ \frac{1}{s}\sup_{t>1}t^{p+2}\exp\left(-a_{\tau,\sigma}g(t)\right). \label{eq.firstBoundsMomentsAbovedc}
\end{align}
We use now an inequality essentially given in~\cite[$(2.3)$]{ptt}:
\begin{equation*}
\sup_{t>1}t^{p}\exp\left(-a_{\tau,\sigma}g(t)\right)
\le A_{1}\sup_{t>1}t^pe^{-T_{\tau,\sigma}(t)}\le A_1\sup_{t>0}t^pe^{-T_{\tau,\sigma}(t)}= A_1M_p^{\tau,\sigma},
\end{equation*}
valid for some $A_1>0$ and every $p\in\N_0$.
With these estimates and property (dc), see~\eqref{eq.Estimate_h_tausigma_dc}, we obtain, continuing with~\eqref{eq.firstBoundsMomentsAbovedc}, that
$$m_{\tau,\sigma}(p)\le \frac{s^{p+1}}{p+1}+ \frac{1}{s}A_1M_{p+2}^{\tau,\sigma}\le \frac{s^{p+1}}{p+1}+ \frac{1}{s}A_1C^{p+2}M_{p}^{\tau,\sigma}.$$
Since $s>0$ was arbitrary, and it is immediate that for any $b>0$ one has
$$
\inf_{s>0}\left(\frac{s^{p+1}}{p+1}+\frac{b}{s}\right)= \frac{p+2}{p+1}b^{(p+1)/(p+2)},
$$
we finally get
$$m_{\tau,\sigma}(p)\le \frac{p+2}{p+1}A_1^{(p+1)/(p+2)} C^{p+1}(M_{p}^{\tau,\sigma})^{(p+1)/(p+2)}\le A_2D_1^{p}M_{p}^{\tau,\sigma}$$
for suitable positive constants $A_2,D_1>0$.

On the other hand, since $g$ is increasing in $(1,\infty)$ by Lemma~\ref{lema.g_hincreasing}, for every $s>0$ we may estimate
$$
m_{\tau,\sigma}(p)\ge \int_0^s t^p\exp\left(-a_{\tau,\sigma}g(t+1)\right)\,dt\ge
\exp\left(-a_{\tau,\sigma}g(s+1)\right)\frac{s^{p+1}}{p+1}.
$$
In particular, we deduce that
$$m_{\tau,\sigma}(p)\ge\sup_{s>1}\exp\left(-a_{\tau,\sigma}g(s+1)\right) \frac{s^{p+1}}{p+1}\ge \frac{1}{2^{2p+1}}\sup_{s>1}(s+1)^{p+1} \exp\left(-a_{\tau,\sigma}g(s+1)\right).$$
We apply now the left inequality in~\eqref{eq.doubleEstimatesSerbian} and obtain
$$
m_{\tau,\sigma}(p)\ge A_{\tau,\sigma}\sup_{s>1}(s+1)^{p+1}e^{-T_{\tau,\sigma}(s+1)}.
$$
Since the sequence $(M_p^{\tau,\sigma})_{p\in\N_0}$ is logarithmically convex, and  $T_{\tau,\sigma}(s)$ is precisely its classical associated function (see~\cite{komatsu}), it is well known that there exists $p_0\in\mathbb{N}$ such that for every $p\ge p_0$ the function $s\in(0,\infty)\mapsto s^{p+1}e^{-T_{\tau,\sigma}(s)}$ reaches its supremum at a point greater than $2$, and the value of such supremum is $M_{p+1}^{\tau,\sigma}$.
Hence, we may deduce the existence of positive constants $E_1,F_1$ such that for every $p\in\mathbb{N}$ one has $m_{\tau,\sigma}(p)\ge E_1F_1^{p}M_{p}^{\tau,\sigma}$.
\end{proof1}

\begin{rem}\label{rem.ChangeBoundsSigmaGreater2}
A version of the previous result can be obtained also for $\sigma\ge 2$, by following an analogous reasoning, but there occur some changes in the form of the estimates. Besides the loss of property (dc), in this case the left-hand side inequality in~\eqref{eq.doubleEstimatesSerbian} suffers from a scaling in the constant $a_{\tau,\sigma}$ appearing in the exponent. As a consequence, we can only prove that there exist $\tilde{K}_1,\tilde{K}_2>0$ and $\tilde{\tau}>\tau$ such that $\tilde{K}_{1}^{p^{\sigma}}M_{p}^{\tau,\sigma}\le m_{\tau,\sigma}(p)\le \tilde{K}_{2}^{p^{\sigma}}M_p^{\tilde{\tau},\sigma}$, for every $p\in\N_0$.

It is also important to remark that the previous result guarantees that the ultraholomorphic  classes associated with the sequences $(m_{\tau,\sigma}(p))_{p\in\N_0}$ and $(M^{\tau,\sigma}_{p})_{p\in\N_0}$ coincide.
\end{rem}

\begin{prop}\label{prop.eOptimalFlatatInfinity}
Let $1<\sigma<2$, $\tau>0$. For every $T\prec S_{2}$ there exist $C_3,K_3,C_4,K_4>0$ (depending on $\tau,\sigma$) such that
\begin{equation*}
C_3\exp(-T_{\tau,\sigma}(|z|/K_3))\le |e_{\tau,\sigma}(z)|\le C_4\exp(-T_{\tau,\sigma}(|z|/K_4)),\quad z\in T.
\end{equation*}
\end{prop}

\begin{proof1}
From Lemma~\ref{lema.Properties.e_h}, there exist $C_2,K_2>0$ such that for every $z\in T$ one has $|e_{\tau,\sigma}(z)|\le C_2e_{\tau,\sigma}(K_2|z|)$. Regarding the values of $e_{\tau,\sigma}(t)$ for positive real $t$, we can apply Lemma~\ref{lema.g_hincreasing}, Lemma~\ref{lema.doubleEstimatesSerbian} and~\eqref{eq.Estimate_h_tausigma_dc} in order to obtain
constants $\tilde{A}_{\tau,\sigma},C>0$ such that
\begin{align*}
e_{\tau,\sigma}(t)&=t\exp(-ag(t+1))<t\exp(-ag(t))\le \tilde{A}_{\tau,\sigma}t\exp(-T_{\tau,\sigma}(t))\\
&=\tilde{A}_{\tau,\sigma}th_{\tau,\sigma}(1/t)\le
\tilde{A}_{\tau,\sigma}t\frac{C}{t}h_{\tau,\sigma}\big(\frac{C}{t}\big)=
C\tilde{A}_{\tau,\sigma}\exp(-T_{\tau,\sigma}(t/C)).
\end{align*}
The inequality on the right is so proved.

Regarding the inequality on the left, again from Lemma~\ref{lema.Properties.e_h} there exist $C_1,K_1>0$ such that for every $z\in T$ one has
$$
|e_{\tau,\sigma}(z)|\ge C_1e_{\tau,\sigma}(K_1|z|)=C_1K_1|z|\exp(-ag(K_1|z|+1)).
$$
Lemma~\ref{lema.g_hincreasing} guarantees that $\lim_{t\to\infty}g'(t)=0$, and the Mean Value Theorem implies then that there exists $c_0>0$ such that for every $t\ge K_1$ one has $g(t+1)-g(t)\le c_0$. So, for every $z\in T$ with $|z|\ge 1$ we have
\begin{align*}
|e_{\tau,\sigma}(z)|&\ge C_1K_1\exp(-ag(K_1|z|))\exp\big(-a(g(K_1|z|+1)-g(K_1|z|))\big)\\
&\ge C_1K_1e^{-ac_0}\exp(-ag(K_1|z|))\ge
C_1K_1e^{-ac_0}A_{\tau,\sigma}\exp(-T_{\tau,\sigma}(K_1|z|)),
\end{align*}
where in the last inequality we have used again Lemma~\ref{lema.doubleEstimatesSerbian}. By suitably enlarging the constants we can make the inequality valid for every $z\in T$.
\end{proof1}

According to~\cite[Theorem~2.3.1]{Thilliez03}, the estimates in Proposition~\ref{prop.eOptimalFlatatInfinity} amount to the fact that the function $e_{\tau,\sigma}(1/z)$ is an optimal flat function in the ultraholomorphic class of functions with a uniform asymptotic expansion associated with the sequence $\M^{\tau,\sigma}$ in any unbounded proper subsector of $\C\setminus(-\infty,0]$.

\section{Right inverses for the asymptotic Borel map in ultraholomorphic classes in sectors}\label{sectRightInver1var}

In this section we construct right inverses for the asymptotic Borel map in ultraholomorphic classes associated with the sequence $\M^{\tau,\sigma}$ in sectors of arbitrary opening of the complex plane.

\begin{theo}\label{tpral}
Let  $\tau>0$ and $1<\sigma<2$ be given. For every $\hat{f}=\sum_{p=0}^{\infty}c_p z^p\in \C[[z]]_{\{\M^{\tau,\sigma}\}}$  there exists a function $f$, holomorphic in $\C\setminus(-\infty,0]$, such that for every unbounded proper subsector $S$ of $\C\setminus(-\infty,0]$, one has that $f$ admits $\hat{f}$ as its uniform asymptotic expansion in $S$.
\end{theo}

\begin{proof1} Let $(m_{\tau,\sigma}(p))_{p\in\N_0}$ be the sequence of moments associated  with the function $e_{\tau,\sigma}(z)$.
Since $\hat{f} \in \C[[z]]_{\{\M^{\tau,\sigma}\}}$, there exist positive constants $C_1,D_1$ such that
$$|c_p|\le C_1D_1^{p} M^{\tau,\sigma}_p,\quad p\in\N_0.$$
From~\eqref{eq.boundsMmomentsdc}, we deduce that the series
$$
\hat{g}=\sum_{p=0}^{\infty}\frac{c_p}{m_{\tau,\sigma}(p)}z^p
$$
is convergent in a disc $D(0,R)$ for some $R>0$, and it defines a holomorphic function $g$ there. Let $0<R_0<R$.
We define
\begin{equation}\label{intope}
f(z):=\int_{0}^{R_0}e_{\tau,\sigma}\left(\frac{u}{z}\right)g(u)\frac{du}{u},\qquad z\in \C\setminus(-\infty,0],
\end{equation}
which may be called a truncated Laplace-like transform of the function $g$ with kernel $e_{\tau,\sigma}$.
By virtue of Leibniz's theorem on analyticity of parametric integrals and the definition of $e_{\tau,\sigma}$, $f$ turns out to be a holomorphic function in $\C\setminus(-\infty,0]$. In order to obtain our result, it suffices to prove that
$f\sim_{\{\M^{\tau,\sigma}\}}\hat{f}$ uniformly in $S_{\delta}$, for every $0<\delta<2$.

Let $N\in\N$ and $z\in S_{\delta}$. We have 
\begin{align*}
f(z)-\sum_{p=0}^{N-1} c_p z^p &= f(z)-\sum_{p=0}^{N-1}\frac{c_p}{m_{\tau,\sigma}(p)} m_{\tau,\sigma}(p)z^p\\
&= \int_{0}^{R_0}e_{\tau,\sigma}\left(\frac{u}{z}\right) \sum_{k=0}^{\infty}\frac{c_{k}}{m_{\tau,\sigma}(k)} u^k \frac{du}{u} -\sum_{p=0}^{N-1}\frac{c_p}{m_{\tau,\sigma}(p)}\int_{0}^{\infty}u^{p-1}e_{\tau,\sigma}(u)du z^p.
\end{align*}
After a change of variable $v=zu$ in the second integral, by virtue of the estimate~(\ref{eq.Bounds_e_h_subsectors_real_axis}) one may use Cauchy's residue theorem in order to check that
$$
z^p\int_{0}^{\infty}u^{p-1}e_{\tau,\sigma}(u)du= \int_{0}^{\infty}v^{p}e_{\tau,\sigma}\left(\frac{v}{z}\right)\frac{dv}{v},
$$
which allows us to write the preceding difference as 
\begin{multline*}
\int_{0}^{R_0}e_{\tau,\sigma}\left(\frac{u}{z}\right) \sum_{k=0}^{\infty}\frac{c_{k}}{m_{\tau,\sigma}(k)} u^k \frac{du}{u} -\sum_{p=0}^{N-1}\frac{c_p}{m_{\tau,\sigma}(p)} \int_{0}^{\infty}u^{p}e_{\tau,\sigma}\left(\frac{u}{z}\right) \frac{du}{u} \\
=\int_{0}^{R_0}e_{\tau,\sigma}\left(\frac{u}{z}\right) \sum_{k=N}^{\infty}\frac{c_{k}}{m_{\tau,\sigma}(k)} u^k\frac{du}{u} -\int_{R_0}^{\infty}e_{\tau,\sigma}\left(\frac{u}{z}\right) \sum_{p=0}^{N-1}\frac{c_p}{m_{\tau,\sigma}(p)} u^{p}\frac{du}{u}.
\end{multline*} 
Then, we have 
$$\left|f(z)-\sum_{p=0}^{N-1} c_p z^p\right|\le f_{1}(z)+f_2(z),$$
where 
$$f_{1}(z)=\left|\int_{0}^{R_0}e_{\tau,\sigma}\left(\frac{u}{z}\right) \sum_{k=N}^{\infty}\frac{c_{k}}{m_{\tau,\sigma}(k)} u^k \frac{du}{u}\right|,$$
$$f_{2}(z)=\left|\int_{R_0}^{\infty}e_{\tau,\sigma}\left(\frac{u}{z}\right) \sum_{p=0}^{N-1}\frac{c_p}{m_{\tau,\sigma}(p)} u^{p}\frac{du}{u}\right|.$$
We now give suitable estimates for $f_1(z)$ and $f_2(z)$.
From~\eqref{eq.boundsMmomentsdc} there exist $C_2,D_2>0$ (not depending on $z$) such that 
\begin{equation}\label{e327}
\frac{c_{k}}{m_{\tau,\sigma}(k)}\le \frac{C_1D_1^{k} M^{\tau,\sigma}_k}{m_{\tau,\sigma}(k)}\le C_2D_2^k,
\end{equation}
for all $k\in\N_0$. This yields
$$f_{1}(z)\le C_2\int_{0}^{R_0}\left|e_{\tau,\sigma}\left(\frac{u}{z}\right)\right| \sum_{k=N}^{\infty}(D_2u)^{k}\frac{du}{u}.$$
Taking $R_0\le(1-\epsilon)/D_2$ for some $\epsilon>0$ if necessary, we get
$$f_1(z)\le\frac{C_2}{\epsilon}D_2^{N}\int_{0}^{R_0} \left|e_{\tau,\sigma}\left(\frac{u}{z}\right)\right|u^{N-1}du.$$

Let us turn our attention to $f_{2}(z)$. We have $u^p\le R_0^pu^N/R_0^N$ for $u\ge R_0$ and $0\le p\le N-1$. So, according to~\eqref{e327}, we may write 
\begin{equation} \label{e1211}
\sum_{p=0}^{N-1}\frac{c_pu^p}{m_{\tau,\sigma}(p)} \leq 
\sum_{p=0}^{N-1} \frac{C_1D_1^pM^{\tau,\sigma}_pu^p}{m_{\tau,\sigma}(p)} 
\leq C_2 \sum_{p=0}^{N-1}D_2^pu^p\le\frac{u^N}{R_0^N} C_2\sum_{p=0}^{N-1}D_2^pR_0^p\le C_5D_5^Nu^N,
\end{equation}
for some positive constants $C_5,D_5$. Then, we conclude
$$f_2(z)\le C_5D^N_5\int_{R_0}^{\infty}\left|e_{\tau,\sigma}\left(\frac{u}{z}\right) \right|u^{N-1}du.$$
So, both $f_1$ and $f_2$ can be estimated above essentially by
$$
\int_{0}^{\infty}\left|e_{\tau,\sigma}\left(\frac{u}{z}\right) \right|u^{N-1}du.
$$
The second inequality in~\eqref{eq.Bounds_e_h_subsectors_real_axis}, a simple change of variable and the estimates in~\eqref{eq.boundsMmomentsdc} imply that, for suitable constants $C_6,K_6,B>0$, one has
\begin{align*}
\int_{0}^{\infty}\left|e_{\tau,\sigma}\left(\frac{u}{z}\right) \right|u^{N-1}du
&\le
\int_{0}^{\infty}C_6e_{\tau,\sigma}\left(K_6\frac{u}{|z|}\right)u^{N-1}du
=C_6K_6^{-N}m_{\tau,\sigma}(N)|z|^N\\
&\le C_6(B/K_6)^{N}M_N^{\tau,\sigma}|z|^N,
\end{align*}
as desired.
\end{proof1}

As we have mentioned in Subsection~\ref{subsec.UltraHolClasses}, for each $A>0$ we know that $(\widetilde{\mathcal{A}}^u_{\{\M\},A}(S), \left\|\,\,\, \right\|_{\M,A,\overset{\sim}{u}})$ and  $(\C[[z]]_{\{\M\},A},\left| \  \right|_{\M,A})$ are Banach spaces. Then, there is a topological consequence of the previous result: for  each $A>0$ we have a  linear and continuous right inverse for the asymptotic Borel map from $\C[[z]]_{\{\M\},A}$ to $\tilde{\mathcal{A}}^u_{\{\M\},dA}(S)$ for some $d>0$.  More precisely, we observe from the proof of Theorem~\ref{tpral} that $D_2=d_1D_1$ for some $d_1\ge1$, in (\ref{e327}). In addition to this, the choice of $R_0$ guarantees that $D_5$ in (\ref{e1211}) is of the form $d_2D_1$ for some $d_2\ge1$. Therefore, one can choose $d:=\frac{2B}{K_6}\max\{d_1,d_2\}$. This entails that if $\hat{f}\in \C[[z]]_{\{\M\},A} $ there exist constants $c,d>0$, not depending neither on $\hat{f}$ nor on $A$, such that for every $N\in\N_0$ one has
$$\left|f(z)-\sum_{p=0}^{N-1} c_p z^p \right| \leq (cC_1)(d A)^N M^{\tau,\sigma}_{N}=
c|\hat{f}|_{\M,A}(dA)^NM^{\tau,\sigma}_N,\qquad z\in S.$$
Then  $\left\| f \right\|_{\M^{\tau,\sigma},dA,\overset{\sim}{u}} \leq c |\hat{f}|_{\M^{\tau,\sigma},A}$, where
the scaling factor $d$ only depends on the unbounded subsector $S$ of $\C\backslash(-\infty,0]$,  that is, $d$ does not depend on $A$.

Consequently, we have the following result which is analogous to a result by V. Thilliez \cite[Thm. 3.2.1]{Thilliez03} for classes associated with strongly regular sequences.

\begin{theo}\label{corolpral}
Let  $\tau>0$ and $1<\sigma<2$ be given. For each $\delta\in (0,2)$ there exists a positive constant $d\ge1$ such that for any $A>0$, the integral operator
$$ T_{\M,A}: \C[[z]]_{\{\M\},A} \longrightarrow \widetilde{\mathcal{A}}^u_{\{\M\},dA}(S_\delta),$$
defined for $\hat{f}= \sum_{p=0}^{\infty}c_p z^p $ in \eqref{intope} by
$$T_{\M,A}(\hat{f}):=\int_{0}^{R_0}e_{\tau,\sigma}(u/z) \Big(\sum_{p=0}^{\infty}\frac{c_{p}}{m_{\tau,\sigma}(p)} u^p\Big) \frac{du}{u}$$
is a linear and continuous right inverse for the asymptotic Borel map $\mathcal{B}$.
\end{theo}

\vskip.3cm
\noindent Affiliations:\par
\noindent ${}^*$ Javier Jim\'enez-Garrido (corresponding author)\\
 Departamento de Matem\'aticas, Estad{\'\i}stica y Computaci\'on,\\
Universidad de Cantabria,\\
Avda. Los Castros, s/n, 39005, Santander, Spain.\\
             Instituto de Investigaci\'on en Matem\'aticas IMUVA, Universidad de Va\-lla\-do\-lid\\
             ORCID: 0000-0003-3579-486X\\
             \email{jesusjavier.jimenez@unican.es}\\
            \vskip.1cm
\noindent Alberto Lastra \\
 Universidad de Alcal\'a\\
Departamento de F\'isica y Matem\'aticas\\
Alcal\'a de Henares, Madrid, Spain. \\
ORCID: 0000-0002-4012-6471\\
\email{alberto.lastra@uah.es}\\
\vskip.1cm
\noindent  Javier Sanz \\
            Departamento de \'Algebra, An\'alisis Matem\'atico, Geometr{\'\i}a y Topolog{\'\i}a\\
            Universidad de Va\-lla\-do\-lid\\
            Facultad de Ciencias, Paseo de Bel\'en 7, 47011 Valladolid, Spain\\
            Instituto de Investigaci\'on en Matem\'aticas IMUVA, Universidad de Valladolid\\
            ORCID: 0000-0001-7338-4971\\
            \email{javier.sanz.gil@uva.es}

\end{document}